\documentclass[graybox]{svmult}

\usepackage{type1cm}        
\usepackage{makeidx}         
\usepackage{graphicx}        
\usepackage{multicol}        
\usepackage[bottom]{footmisc}

\usepackage{newtxtext}       %
\usepackage{newtxmath}       

\usepackage{cite}
\usepackage{multirow}

\usepackage[ruled, lined, longend, linesnumbered]{algorithm2e}

\makeindex             

\begin{document}

\title*{Model reduction by separation of variables: a comparison between Hierarchical Model reduction and Proper Generalized Decomposition}
\titlerunning{Model reduction by separation of variables: a comparison between HiMod and PGD}
\author{Simona Perotto, Michele Giuliano Carlino, Francesco Ballarin}
\institute{
Simona Perotto \at MOX - Modeling and Scientific Computing, Dipartimento di Matematica, Politecnico di Milano, Piazza Leonardo da Vinci 32, I-20133 Milano, Italy \email{simona.perotto@polimi.it}
\and
Michele Giuliano Carlino \at Inria Bordeaux Sud-Ouest and Institut de Math\'ematiques de Bordeaux, University of Bordeaux,
351, cours de la Lib\'eration, F-33405 Talence \email{michele-giuliano.carlino@inria.fr}
\and
Francesco Ballarin \at mathLab, Mathematics Area, SISSA, via Bonomea 265, I-34136 Trieste, Italy \email{francesco.ballarin@sissa.it}}

\maketitle

\abstract{Hierarchical Model reduction and Proper Generalized 
Decomposition both exploit separation of variables
to perform a model reduction. After setting the basics, we exemplify 
these techniques on some standard elliptic problems
to highlight pros and cons of the two procedures, both from a 
methodological and a numerical viewpoint.}

\section{Introduction}
This paper is meant as a first attempt to compare two procedures which share the idea of exploiting separation of 
variables to perform model reduction, albeit with different purposes. 
Proper Generalized Decomposition (PGD) is essentially employed as a powerful tool to deal with parametric 
problems in several fields of application~\cite{AmmarCuetoChinesta12,NiroomandietAl13,SignoriniZlotnikDiez17}.
Parametrized models characterize multi-query contexts, such as parameter optimization,
statistical analysis or inverse problems. Here, the computation of the solution for many different parameters
demands, in general, a huge computational effort, and this justifies the development of model reduction techniques.

For this purpose, projection-based techniques, such as Proper Orthogonal Decomposition (POD) or Reduced Basis
methods, are widely used in the literature~\cite{HesthavenRozzaStamm16}. The idea is to project the discrete operators onto a 
reduced space so that the problem can be solved rapidly in the lower dimensional space.
PGD adopts a completely different way to deal with parameters. Here, parameters are considered as new independent
variables of the problem, together with the standard space-time ones~\cite{ChinestaKeuningsLeygue13}. 
Although the dimensionality of the problem is inevitably increased, PGD transforms the computation of the solution for new values of the parameters into a plain evaluation of the reduced solution, with striking computational advantages.\\
Hierarchical-Model (HiMod) reduction has been proposed to improve one-di\-men\-sion\-al (1D) partial differential equation (PDE) solvers for problems defined in domains with a geometrically dominant direction, like slabs or pipes~\cite{ErnPerottoVeneziani08,PerottoErnVeneziani10}.
The main applicative field of interest is hemodynamics, in particular the modeling of blood flow in patient-specific geometries.
Purely 1D hemodynamic models completely drop the transverse dynamics, which, however may be locally important (e.g.,
in the presence of a stenosis or an aneurism). HiMod aims at providing a numerical tool to incorporate the transverse components of the 3D solution into a conceptually 1D solver. To do this, the driving idea is to discretize main and transverse dynamics in a different way. The latter are generally of secondary importance and can be described by few degrees of freedom using a spectral approximation, in combination, for instance, with a finite element (FE) discretization of the mainstream.

The parametric version of HiMod (namely, HiPOD) 
is a more recent proposal~\cite{LupoPasiniPerottoVeneziani19,BarolietAl17}. On the other hand, PGD is not so widely employed 
in a non-parametric setting, despite its original formulation~\cite{LadevezePassieuxNeron10}. 
Nevertheless, for the sake of comparison, in this paper we consider the non-parametric as well as the parametric versions of 
both the HiMod and PGD approaches.
The goal is to begin a preliminary comparative analysis between the two 
methodologies, to highlight the respective weaknesses and strengths.
The main limit of PGD remains its inability to deal with non-Cartesian geometries without losing the computational benefits arising from the separability of the spatial coordinates. HiMod turns out to be more flexible from a geometric viewpoint. On the other hand, PGD turns out to be extremely effective for parametric problems thanks to the explicit expression of the PGD solution in terms of the parameters, while HiPOD can be classified as a projection-based method with all the associated drawbacks.
In perspective, the ultimate goal is to merge HiMod with PGD to emphasize the good features and mitigate the intrinsic limits of the two methods taken alone.

\section{The HiMod approach}\label{sec:himod}
Hierarchical Model reduction proved to be an efficient and reliable method to deal with phenomena 
charaterized by dominant dynamics~\cite{GuzzettiPerottoVeneziani18}. In general,
the computational domain itself exhibits an intrinsic directionality. We assume 
$\Omega \subset \mathbb{R}^d$ ($d = 2, 3$) to coincide with a $d$-dimensional fiber bundle,
$\Omega = \bigcup_{x \in \Omega_{1D}} \{ x \} \times \gamma_x$, where $\Omega_{1D}\subset \mathbb{R}$ 
denotes the supporting fiber aligned with the main stream, while $\gamma_x \subset \mathbb{R}^{d-1}$ is the 
transverse fiber at $x\in \Omega_{1D}$, parallel to the transverse dynamics.
For the sake of simplicity, we identify $\Omega_{1D}$ with a straight segment, $(x_0, x_1)$. 
We refer to~\cite{Perotto14,PerottoRealiRusconiVeneziani17} for the case where $\Omega_{1D}$ is curvilinear.
From a computational viewpoint, the idea is to exploit a map, $\Psi:\Omega \rightarrow \hat \Omega$,
transforming the physical domain, $\Omega$, into a reference domain, $\hat \Omega$, 
and to make explicit computations in $\hat \Omega$ only. 
Typically,  $\hat \Omega$ coincides with a rectangle in 2D, with a cylinder with circular section in 3D.
To define $\Psi$, for each $x\in \Omega_{1D}$, we introduce the map, $\psi_x: \gamma_x \rightarrow \hat{\gamma}_{d-1}$,
from fiber $\gamma_x$ to the reference transverse fiber, $\hat{\gamma}_{d-1}$, so that 
the reference domain coincides with $\hat{\Omega} = \bigcup_{x \in \Omega_{1D}} \{ x \} \times \hat{\gamma}_{d-1}$. The supporting fiber is preserved by map $\Psi$, which modifies the lateral boundaries only.

We consider now the (full) problem to be reduced. Due to the comparative purposes of the paper, we 
focus on a scalar elliptic equation, and, in particular, on the associated weak formulation,
\begin{equation}\label{eq:fullweak}
	\text{find } u \in V: a(u,v) = F(v) \quad \forall v \in V,
\end{equation}
where $V \subseteq H^1(\Omega)$, $a(\cdot, \cdot): V \times V \rightarrow \mathbb{R}$ is a continuous and coercive bilinear form and $F(\cdot): V \rightarrow \mathbb{R}$ is a continuous linear functional.
To provide the HiMod formulation for problem \eqref{eq:fullweak}, we introduce the hierarchical reduced space
\begin{equation}\label{eq:himodspace}
	V_m = \Bigg\{ v_m(x, \vec{y}) = \sum_{k = 1}^m \tilde{v}_k(x) \varphi_k(\psi_x(\vec{y})), \text{ with } \tilde{v}_k \in V_{1D}^h, \, x \in \Omega_{1D}, \, \vec{y} \in \gamma_x \Bigg\}
\end{equation}
for a modal index $m\in \mathbb{N}^+$, 
where $V_{1D}^h \subseteq H^1(\Omega_{1D})$ is a discrete space of dimension $N_h$
associated with a partition $\mathcal{T}_h$ of $\Omega_{1D}$, 
while $\{ \varphi_k \}_{k=1}^m$ denotes a modal basis of functions orthogonal with respect to the $L^2(\hat{\gamma}_{d-1})$-scalar product. Index $m$ sets the hierarchical level of the HiMod space, being $V_m \subset V_{m+1}$, for any $m$. 
Concerning $V_{1D}^h$, we adopt here a standard FE space, although any discrete space can be employed (see, e.g.,~\cite{PerottoRealiRusconiVeneziani17}, where an isogeometric discretization is used). Functions in $V_{1D}^h$ have to include the boundary conditions on $\{x_0\}\times \gamma_{x_0}$ and $\{x_1\}\times \gamma_{x_1}$; analogously, the modal functions have to take into account the boundary data along the horizontal sides. 
In Sect.~\ref{sec4} further comments are provided about the selection of the modal basis and of the modal index $m$.
The HiMod formulation for problem \eqref{eq:fullweak} thus reads
\begin{equation}\label{eq:himodweak}
	\text{find } u_m^{\rm HiMod} \in V_m: a(u_m^{\rm HiMod},v_m) = F(v_m) \quad \forall v_m \in V_m.
\end{equation} 
To ensure the well-posedness of formulation (\ref{eq:himodweak}) and
the convergence of the HiMod approximation, $u_m^{\rm HiMod}$, to the full solution, $u$, we 
endow the HiMod space with a conformity and a spectral approximability hypothesis, and we introduce 
a standard density assumption on the discrete space $V_{1D}^h$ (see~\cite{PerottoErnVeneziani10} for all the details).

The HiMod solution can be fully characterized by introducing 
a basis, $\{ \theta_l \}_{l=1}^{N_h}$, for the space $V_{1D}^h$. Actually, each modal coefficient,
$\tilde{u}_k$, of $u_m^{\rm HiMod}$ can be expanded in terms of such a basis, so that, we obtain the modal representation
\begin{equation}\label{eq:himodsolution}
	u_m^{\rm HiMod}(x, \vec{y}) = \sum_{k=1}^m \sum_{l=1}^{N_h} \tilde{u}_{k,l} \theta_l(x) \varphi_k(\psi_x(\vec{y})).
\end{equation}
The actual unknowns of problem \eqref{eq:himodweak} become the $mN_h$ coefficients $\{ \tilde{u}_{k,l} \}_{k=1, l=1}^{m, N_h}$. With reference to the Poisson problem, $-\Delta u=f$,
completed with full homogeneous Dirichlet boundary data, the corresponding HiMod formulation,
after exploiting \eqref{eq:himodsolution} in \eqref{eq:himodweak} and picking $v_m(x, \vec{y})=\theta_i(x) \varphi_j(\psi_x(\vec{y}))$ with $i=1, \ldots, N_h$ and $j=1, \ldots, m$,
reduces to the system of $mN_h$ 1D equations in the $mN_h$ unknowns $\{ \tilde{u}_{k,l} \}_{k=1, l=1}^{m, N_h}$,
$$
\begin{array}{l}
\displaystyle  \sum_{k=1}^m \sum_{l=1}^{N_h} \tilde{u}_{k,l} \Bigg[ \displaystyle \int_{\Omega_{1D}} \bigg( \hat r_{jk}^{1, 1}(x)
\frac{d\theta_l}{dx}(x)\frac{d\theta_i}{dx}(x) + \hat r_{jk}^{1, 0}(x)\frac{d\theta_l}{dx}(x) \theta_i(x)+ \\[4mm]
\hspace*{1.5cm} + \, \hat r_{jk}^{0, 1}(x) \theta_l(x) \displaystyle \frac{d\theta_i}{dx}(x) +  \hat r_{jk}^{0, 0}(x)  \theta_l(x) \theta_i(x) \bigg) \Bigg] \, dx =  \int_{\Omega_{1D}} \hat f_j (x) \theta_i(x) \, dx,
 \end{array}
$$
where $\hat r^{a, b}_{jk}(x)=\int_{\hat \gamma_{d-1} }r^{a, b}_{jk}(x, \hat{\vec{y}})|J|\, d\hat{\vec{y}}$ with $a, b=0$,$1$, $J={\rm det}\big( {\mathcal D}_2^{-1}(x, \psi_x^{-1}(\hat{\vec{y}}) )\big)$ with ${\mathcal D}_2={\mathcal D}_2(x, \psi_x^{-1}(\hat{\vec{y}}) )=\nabla_{\vec{y}}\psi_x$,
$$
\begin{array}{lll}
r^{0, 0}_{jk}(x, \hat{\vec{y}})=\varphi'_k(\hat{\vec{y}})\varphi'_j(\hat{\vec{y}}) \big({\mathcal D}_1^2 + {\mathcal D}_2^2\big),\quad 
&r^{0, 1}_{jk}(x, \hat{\vec{y}})=\varphi'_k(\hat{\vec{y}})\varphi_j(\hat{\vec{y}}) {\mathcal D}_1,
\\[2mm] 
r^{1, 0}_{jk}(x, \hat{\vec{y}})=\varphi_k(\hat{\vec{y}})\varphi'_j(\hat{\vec{y}}) {\mathcal D}_1,\quad
&r^{1, 1}_{jk}(x, \hat{\vec{y}})=\varphi_k(\hat{\vec{y}})\varphi_j(\hat{\vec{y}}),
\end{array}
$$
with ${\mathcal D}_1={\mathcal D}_1(x, \psi_x^{-1}(\hat{\vec{y}}) )=\partial \psi_x/\partial x$,
and $\hat f_j (x)=\int_{\hat \gamma_{d-1} } f(x, \psi_x^{-1}(\hat{\vec{y}}) ) \varphi_j(\hat{\vec{y}}) |J|\, d \hat{\vec{y}}$. Information associated with the transverse dynamics are lumped in the coefficients 
$\{ \hat r^{a, b}_{jk} \}$, so that the HiMod system is solved on the supporting fiber, $\Omega_{1D}$.
Collecting the HiMod unknowns, by mode, in the vector $\vec{u}_m^{\rm HiMod}\in \mathbb{R}^{m N_h}$, such that
\begin{equation}\label{mod_expansion}
\vec{u}_m^{\rm HiMod}=[\tilde{u}_{1,1}, \tilde{u}_{1,2}, \ldots, \tilde{u}_{1,N_h}, \tilde{u}_{2,1}, \ldots, 
\tilde{u}_{m,1}, \ldots, \tilde{u}_{m,N_h}]^T,
\end{equation} 
we can rewrite the HiMod system in the compact form
\begin{equation}\label{eq:himodalgebraic}
	A_m^{\rm HiMod}  \vec{u}_m^{\rm HiMod} = \vec{f}_m^{\rm HiMod},
\end{equation}
where $A_m^{\rm HiMod} \in \mathbb{R}^{m N_h \times m N_h}$ and $\vec{f}_m^{\rm HiMod} \in \mathbb{R}^{m N_h}$ are the HiMod stiffness matrix and right-hand side, respectively, with 
$[ \vec{f}_m^{\rm HiMod} ]_{ji}=\int_{\Omega_{1D}} \hat f_j (x) \theta_i(x) dx$, and
$[A_m^{\rm HiMod}]_{ji,kl}=
\sum_{a,b=0}^1 \int_{\Omega_{1D}} \hat r^{a, b}_{jk}(x) \frac{d^a \theta_l}{dx}(x) 
\frac{d^b \theta_i}{dx}(x) dx$. 
According to \eqref{mod_expansion}, for each modal index $j$, between $1$ and $m$, the nodal index, $i$, takes the values $1, \ldots, N_h$.
Thus, HiMod reduction leads to solve a system of order $mN_h$, independently of the dimension of the full problem \eqref{eq:fullweak}.

\section{The PGD approach}\label{sec:pgd}
To perform PGD, we have to introduce on problem \eqref{eq:fullweak}
a separability hypothesis with respect to both the spatial variables and the 
data~\cite{ChinestaKeuningsLeygue13,PruliereChinestaAmmar10}. 
Thus, domain $\Omega\subset \mathbb{R}^d$ coincides with
the rectangle $\Omega_x \times \Omega_y$ if $d = 2$, with
the parallelepiped $\Omega_x \times \Omega_y \times \Omega_z$ (total separability) 
or with the cylinder $\Omega_x \times \Omega_{\vec{y}}$ (partial separability) if $d=3$, 
for $\Omega_x$, $\Omega_y$, $\Omega_z \subset \mathbb{R}$ and  
$\Omega_{\vec{y}} \subset \mathbb{R}^{2}$, being ${\bf y}=(y, z)$.
In the following, we focus on partial separability, since it is more suited to match 
HiMod reduction with PGD. Analogously, we assume that the generic problem data, $d=d(x, y, z)$,  can be
written as $d=d^x(x)d^{\bf y}({\bf y})$.
The separability is inherited by the PGD space
\begin{equation}\label{eq:pgdspace}
W_m = \Bigg\{ w_m(x, \vec{y}) = \sum_{k = 1}^m w_k^x(x) w_k^{\vec{y}}(\vec{y}), \text{ with }  w_k^x \in W^x_h, \, w_k^{\vec{y}} \in 
W^{\vec{y}}_h, \, x \in \Omega_x, \, \vec{y} \in \Omega_{\vec{y}} \Bigg\},
\end{equation}
where $W^x_h \subseteq H^1(\Omega_x)$ and 
$W^{\vec{y}}_h \subseteq H^1(\Omega_{\vec{y}}; \mathbb{R}^{d-1})$
are discrete spaces, with $\text{dim}(W^x_h)=N_h^x$ and $\text{dim}(W^{\vec{y}}_h)=N_h^{\vec{y}}$, associated with 
partitions, $\mathcal{T}_h^x$ and $\mathcal{T}_h^{\vec{y}}$, of $\Omega_x$ and $\Omega_{\vec{y}}$, respectively. 
In general, $W^x_h$ and $W^{\vec{y}}_h$ are FE spaces, although, a priori, any discretization 
can be adopted. It turns out that $W_m$ is a tensor function space, being 
$W_m = W^x_h \otimes W^{\vec{y}}_h \subseteq H^1(\Omega_x) \otimes H^1(\Omega_{\vec{y}}; \mathbb{R}^{d-1})$.\\ 
Index $m$ plays the same role as in the HiMod reduction, setting the level of detail for the reduced solution
(see Sect.~\ref{sec4} for possible criteria to choose $m$).
PGD exploits the hierarchical structure in $W_m$ to build the generic function $w_m\in W_m$.
In particular, $w_m$ is computed as
\begin{equation}\label{pgd_rep}
w_m(x, \vec{y}) = w_m^x(x) w_m^{\vec{y}}(\vec{y}) + \sum_{k = 1}^{m-1} w_k^x(x) w_k^{\vec{y}}(\vec{y}),
\end{equation}
where $w_k^x$ and $w_k^{\vec{y}}$ are assumed known for $k=1, \ldots, m-1$, so that  
the enrichment functions, $w_m^x$ and $w_m^{\vec{y}}$, become the actual unknowns. 
To provide the PGD formulation for the Poisson problem considered in Sect.~\ref{sec:himod}, 
we exploit representation \eqref{pgd_rep} for the PGD approximation, $u^{\rm PGD}_m$, and we pick the test function
as $X(x) Y(\vec{y})$,
with $X\in W_h^x$ and $Y\in W^{\vec{y}}_h$. The coupling between the unknowns, $u_m^x$ and $u_m^{\vec{y}}$, leads to a nonlinear problem, which is tackled by means of the \textit{Alternating Direction Strategy} (ADS)~\cite{ChinestaKeuningsLeygue13}.
The idea is to look for $u_m^x$ and $u_m^{\vec{y}}$, separately via a fixed point procedure.
We introduce an auxiliary index to keep trace of the ADS iterations, so that,
at the $p$-th ADS iteration we compute $u_m^{x, p}$ and $u^{\vec{y},p}_m$ starting from the previous
approximations, $u_m^{x, g}$ and $u^{\vec{y},g}_m$ for $g=1, \ldots, p-1$, following a two-step procedure.
First, we compute $u_m^{x, p}$ by identifying $u_m^{\vec{y}}$ with $u^{\vec{y},p-1}_m$, 
and by selecting $Y(\vec{y})=u^{\vec{y},p-1}_m$ in the test function. This yields, for any $X\in W_h^x$, 
\begin{equation}\label{system1}
\begin{array}{ll}
&\int_{\Omega_x} \big( u_m^{x, p} \big)' X'dx\int_{\Omega_{\vec{y}}}\big[ u^{\vec{y},p-1}_m\big]^2 d{\vec{y}}+
\int_{\Omega_x} u_m^{x, p} X dx\int_{\Omega_{\vec{y}}}\big[ \big( u^{\vec{y},p-1}_m \big)'\big]^2 d{\vec{y}}\\[2mm]
=& \int_{\Omega_x} f^xXdx \int_{\Omega_{\vec{y}}} f^{\vec{y}} u^{\vec{y},p-1}_md{\vec{y}} - \sum_{k=1}^{m-1}
\int_{\Omega_x} \big( u_k^{x} \big)' X'dx\int_{\Omega_{\vec{y}}} u_k^{\vec{y}} u^{\vec{y},p-1}_m d{\vec{y}}\\[2mm]
-& \sum_{k=1}^{m-1}
\int_{\Omega_x} u_k^{x} X dx\int_{\Omega_{\vec{y}}} \big( u_k^{\vec{y}} \big)' \big( u^{\vec{y},p-1}_m \big)' d{\vec{y}}, 
\end{array}
\end{equation}
where the separability of $f$ is exploited (the dependence on the independent variables, $x$ and $\vec{y}$, is omitted to simplify notation).
Successively, we compute $u_m^{\vec{y}, p}$, after setting $u_m^x$ to $u_m^{x,p}$ and choosing function $X$  
as to $u_m^{x,p}$ in the test function, so that we obtain, for any $Y\in W_h^{\vec{y}}$,
\begin{equation}\label{system2}
\begin{array}{ll}
&\int_{\Omega_x}\big[ \big( u_m^{x, p} \big)' \big]^2 dx\int_{\Omega_{\vec{y}}}u^{\vec{y},p}_m Y d{\vec{y}}+
\int_{\Omega_x} \big[ u_m^{x, p} \big]^2 dx\int_{\Omega_{\vec{y}}}\big( u^{\vec{y},p}_m \big)' Y' d{\vec{y}}\\[2mm]
=& \int_{\Omega_x} f^x u_m^{x, p} dx \int_{\Omega_{\vec{y}}} f^{\vec{y}} Y d{\vec{y}} - \sum_{k=1}^{m-1}
\int_{\Omega_x} \big( u_k^{x} \big)' \big( u_m^{x, p} \big)'dx\int_{\Omega_{\vec{y}}} u_k^{\vec{y}} Y d{\vec{y}}\\[2mm]
-& \sum_{k=1}^{m-1}
\int_{\Omega_x} u_k^{x} u_m^{x, p} dx\int_{\Omega_{\vec{y}}} \big( u_k^{\vec{y}} \big)' Y' d{\vec{y}}.
\end{array}
\end{equation}
The algebraic counterpart of \eqref{system1} and \eqref{system2} is obtained by introducing a basis, $
{\mathcal B}_x=\{ \theta_\alpha^x \}_{\alpha = 1}^{N_h^x}$ and ${\mathcal B}_{\vec{y}}=\{ \theta_\beta^{\vec{y}} \}_{\beta = 1}^{N_h^{\vec{y}}}$, for the space $W^x_h$ and $W^{\vec{y}}_h$, respectively, so that
$u_j^q(q)=\sum_{i=1}^{N_h^q} \tilde u_{ji}^q \theta_i^q(q)$, $u_m^{q,s}(q)=\sum_{i=1}^{N_h^q} \tilde u_{mi}^{q,s} \theta_i^q(q)$,
with $q=x$, $\vec{y}$, $s=p$, $p-1$, $j=1,\ldots, m-1$, and, likewise, $X(x)=\sum_{\alpha=1}^{N_h^x} \tilde x_{\alpha} \theta_\alpha^x(x)$ 
and $Y(\vec{y})=\sum_{\beta=1}^{N_h^{\vec{y}}} \tilde y_{\beta} \theta_{\beta}^{\vec{y}}(\vec{y})$. Thanks to these expansions and to the 
arbitrariness of $X$ and $Y$, we can rewrite \eqref{system1} and \eqref{system2} as
\begin{equation}\label{sys1}
\begin{array}{ll}
&\Big\{ \Big[ \big( {\bf u}_m^{\vec{y}, p-1}\big)^T M^{\vec{y}} {\bf u}_m^{\vec{y}, p-1}\Big] K^{x} + \Big[ \big( {\bf u}_m^{\vec{y}, p-1}\big)^T K^{\vec{y}} {\bf u}_m^{\vec{y}, p-1}\Big] M^x \Big\} {\bf u}_m^{x, p} = \Big[ \big( {\bf u}_m^{\vec{y}, p-1}\big)^T {\bf f}^{\vec{y}}\Big]  {\bf f}^{x}\\[2mm] 
- &\sum_{k=1}^{m-1}
\Big\{ \Big[ \big( {\bf u}_m^{\vec{y}, p-1}\big)^T M^{\vec{y}} {\bf u}_k^{\vec{y}}\Big] K^{x} + \Big[ \big( {\bf u}_m^{\vec{y}, p-1}\big)^T K^{\vec{y}} {\bf u}_k^{\vec{y}}\Big] M^x \Big\} {\bf u}_k^{x},
\end{array}
\end{equation}
and
\begin{equation}\label{sys2}
\begin{array}{ll}
&\Big\{ \Big[ \big( {\bf u}_m^{x, p}\big)^T K^x {\bf u}_m^{x, p}\Big] M^{\vec{y}}+ \Big[ \big( {\bf u}_m^{x, p}\big)^T M^{x} {\bf u}_m^{x, p}\Big] K^{\vec{y}} \Big\} {\bf u}_m^{\vec{y}, p} = \Big[ \big( {\bf u}_m^{x, p}\big)^T {\bf f}^{x}\Big]  {\bf f}^{\vec{y}}\\[2mm] 
- &\sum_{k=1}^{m-1}
\Big\{ \Big[ \big( {\bf u}_m^{x, p}\big)^T K^{x} {\bf u}_k^{x}\Big] M^{\vec{y}} + \Big[ \big( {\bf u}_m^{x, p}\big)^T M^{x} {\bf u}_k^{x}\Big] K^{\vec{y}} \Big\} {\bf u}_k^{\vec{y}},
\end{array}
\end{equation}
respectively, where vectors ${\bf u}_j^q$, ${\bf u}_m^{q,s}\in \mathbb{R}^{N_h^q}$ 
collect the PGD coefficients, being $\big[ {\bf u}_j^q \big]_i=\tilde u_{ji}^q$, 
$\big[ {\bf u}_m^{q,s} \big]_i=\tilde u_{mi}^{q,s}$ and $i=1,\ldots, N_h^q$, 
$K^x$, $M^x\in \mathbb{R}^{N_h^x\times N_h^x}$ and $K^{\vec{y}}$, 
$M^{\vec{y}}\in \mathbb{R}^{N_h^{\vec{y}}\times N_h^{\vec{y}}}$ are the stiffness and mass matrices 
associated with $x$- and ${\vec{y}}$-variables, with $\big[ K^x \big]_{\alpha l}=\int_{\Omega_x} \big( \theta_\alpha^x \big)' \big( \theta_l^x \big)' dx$,
$\big[ K^{\vec{y}} \big]_{\beta s}=\int_{\Omega_{\vec{y}}} \big( \theta_\beta^{\vec{y}} \big)' \big( \theta_s^{\vec{y}} \big)' d{\vec{y}}$,
$\big[ M^x \big]_{\alpha l}=\int_{\Omega_x}  \theta_\alpha^x \theta_l^x dx$,
$\big[ M^{\vec{y}} \big]_{\beta s}=\int_{\Omega_{\vec{y}}} \theta_\beta^{\vec{y}} \theta_s^{\vec{y}} d{\vec{y}}$, and 
where ${\bf f}^x\in \mathbb{R}^{N_h^x}$, ${\bf f}^{\vec{y}} \in \mathbb{R}^{N_h^{\vec{y}} }$, 
with $\big[ {\bf f}^x\big]_l=\int_{\Omega_x} f^x \theta_l^x dx$, 
$\big[ {\bf f}^{\vec{y}} \big]_s=\int_{\Omega_{\vec{y}}} f^{\vec{y}} \theta_s^{\vec{y}} d{\vec{y}}$, for $\alpha, l=1, \ldots, N_h^x$, $\beta, s=1, \ldots, N_h^{\vec{y}}$.
Systems \eqref{sys1} and \eqref{sys2} are solved at each ADS iteration, so that the computational 
effort characterizing PGD is the one associated with the solution of two systems 
of order $N_h^x$ and $N_h^{\vec{y}}$, respectively, for each ADS iteration.
When a certain stopping criterion is met (see the next section for more details), ADS procedure yields vectors 
${\bf u}_m^x$ and ${\bf u}_m^{\vec{y}}$ which identify the enrichment functions $u_m^x$ and $u_m^{\vec{y}}$.


\section{HiMod reduction versus PGD}\label{sec4} 
Both HiMod reduction and PGD exploit the separation of variables and, according to~\cite{ChinestaKeuningsLeygue13},
belong to the {\it a priori} approaches, since they do not
rely on any solution to the problem at hand. Nevertheless, we can easily itemize 
features which distinguish the two techniques. The most relevant ones concern the geometry of $\Omega$,
the selection of the transverse basis and of the modal index, and the numerical implementation of the two procedures.
Pros and cons of the two methods are then here highlighted. 

\subsection{Domain geometry}
HiMod reduction and PGD advance precise hypotheses on the geometry of the computational domain.

According to the HiMod approach, $\Omega$ is expected to coincide with a fiber bundle and to be mapped into 
the reference domain, $\hat \Omega$, by a sufficiently regular transformation. Actually, map $\Psi$ is assumed 
differentiable, while map $\psi_x$ is required to be a $C^1$-diffeomorphism, for all $x\in \Omega_{1D}$~\cite{PerottoErnVeneziani10}.
These hypotheses introduce some constraints, in particular, on the lateral boundary of 
$\Omega$ which, e.g., cannot exhibit kinks. Additionally, geometries of interest in many applications, 
such as bifurcations or, more in general, networks are ruled out from the demands on $\psi_x$ and $\Psi$. 
An approach based on the domain decomposition technique is currently under investigation as a viable way to deal with 
such geometries. 
The isogeometric version of HiMod (i.e., the HIgaMod approach) will play a crucial role
in view of HiMod simulations for the blood flow modeling in patient-specific geometries~\cite{PerottoRealiRusconiVeneziani17}.

The constraints introduced by PGD on the geometry of $\Omega$ are more restrictive. The separability hypothesis leads to 
consider essentially only Cartesian domains. This considerably reduces the applicability of PGD to practical contexts. 
Some techniques are available in the literature to overcome this issue. For instance, 
in~\cite{GonzalesetAl10} a generic domain is embedded into a Cartesian geometry, while
in~\cite{GhantiosetAl12} the authors introduce a parametrization map for quadrilateral domains.

Overall, HiMod reduction exhibits a higher geometric flexibility with respect to PGD, in its
straightforward formulation. As discussed in Sect.~\ref{parametri_sec}, this limitation can be removed when considering a
parametric setting. 

\subsection{Modeling of the transverse dynamics}
In the HiMod expansion, $\vec{y}$-components, $\varphi_k(\psi_x(\vec{y}))$, are selected 
before starting the model reduction. This choice, although coherent with an a priori approach, 
introduces a constraint on the dynamics that can be described, so that
hints about the solution trend along the transverse direction can be helpful to select a representative modal basis. 
In the original proposal of the HiMod procedure, sinusoidal functions are employed according to a Fourier 
expansion~\cite{ErnPerottoVeneziani08,PerottoErnVeneziani10}. 
This turns out to be a reasonable choice when Dirichlet boundary conditions are assigned on the lateral surface,
$\Gamma_{\rm lat}=\{ x \} \times \partial \gamma_x$, of $\Omega$.  
Legendre polynomials, properly modified to include the homogeneous Dirichlet data and orthonormalized, are 
employed in~\cite{PerottoErnVeneziani10} as an alternative to a trigonometric expansion.
Nevertheless, Legendre polynomials require high-order quadrature rules to accurately compute coefficients $\{ \hat r^{a, b}_{jk} \}$. \\
In~\cite{AlettiPerottoVeneziani18}, the concept of educated modal basis is introduced to impose generic boundary conditions on $\Gamma_{\rm lat}$. 
The idea is to solve an auxiliary Sturm-Liouville eigenvalue problem on the transverse reference fiber $\hat \gamma_{d-1}$, 
to build a basis which automatically includes the boundary values on $\Gamma_{\rm lat}$. The eigenfunctions
of the Sturm-Liouville problem provide the modal basis.
A first attempt to generalize the educated-HiMod reduction to three-dimensional (3D) cylindrical geometries 
is performed in~\cite{GuzzettiPerottoVeneziani18}, where the Navier-Stokes
equations are hierarchically reduced to model the blood flow in pipes. 
This generalization is far from being straightforward
due to the employment of polar coordinates. 
To overcome this issue, we are currently investigating the HIgaMod approach~\cite{PerottoRealiRusconiVeneziani17}, which allows us to define 
the transverse basis as the Cartesian product 
of 1D modal functions, independently of the considered geometry.\\
Additionally, we remark that any modal basis can be precomputed on the transverse reference fiber before performing the HiMod reduction, 
thanks to the employment of map $\Psi$. This considerably simplifies computations. 

When applying PGD, $\vec{y}$-components are unknown as the ones associated with $x$. 
This leads to the nonlinear problems \eqref{system1}-\eqref{system2}, thus loosing any
advantage related to a precomputation of the HiMod modal basis. On the other hand, 
PGD does not constrain the transverse dynamic to follow a prescribed (e.g., sinusoidal) analytical shape
as HiMod procedure does. The educated-Himod reduction clearly is out of this comparison, since the modal basis 
strictly depends on the problem at hand.

Finally, we observe that HiMod modes are orthonormal with respect to the $L^2(\hat \gamma_{d-1})$-norm. 
This property is not ensured by PGD. 

Concerning the selection of the modal index $m$ in \eqref{eq:himodspace} and \eqref{eq:pgdspace}, as a first attempt,
both HiMod reduction and PGD resort to a trial-and-error approach, so that the modal index is gradually increased until
a check on the accuracy of the reduced solution is satisfied. 
For instance, in~\cite{ErnPerottoVeneziani08,PerottoErnVeneziani10} a qualitative 
investigation of the contour plot of the HiMod approximation drives the choice of $m$. 
Concerning PGD, the check on the relative enrichment
\begin{equation}\label{m_criterion}
\frac{\| u_m^x u_m^{\vec{y}} \|_{L^2(\Omega)}}{\| u_1^x u_1^{\vec{y}} \|_{L^2(\Omega)}}\le {\tt TOL_E},
\end{equation}
is usually employed, with ${\tt TOL_E}$ a user-defined tolerance~\cite{ChinestaKeuningsLeygue13}.
An automatic selection of index $m$ can yield a significant improvement. 
In~\cite{PerottoVeneziani14,PerottoZilio15}, an adaptive procedure is proposed for HiMod, based on 
an a posteriori modeling error analysis. In particular, the estimator in~\cite{PerottoVeneziani14} is derived 
in a goal-oriented setting to control a quantity of interest, and exploits the hierarchical structure 
(i.e., the inclusion $V_m\subset V_{m+d}$, $\forall m$, $d\in \mathbb{N}^+$) typical of a HiMod reduction. 
A similar modeling error analysis is performed in~\cite{AmmaretAl10} for PGD, although 
no adaptive algorithm is here set to automatically pick the reduced model.
Paper~\cite{PerottoZilio15} generalizes the a posteriori analysis in~\cite{PerottoVeneziani14} to an unsteady setting, providing the tool
to automatically select $m$ together with the partition $\mathcal{T}_h$ along $\Omega_{1D}$ and the time step.  

Finally, HiMod allows to tune the modal index along the domain $\Omega$, according to the local 
complexity of the transverse dynamics. In particular, $m$ can be varied in different areas of $\Omega$ or, in the presence
of very localized dynamics, in correspondence with specific nodes of the partition $\mathcal{T}_h$.
We refer to these two variants as to {\it piecewise} and {\it pointwise} HiMod reduction, in contrast to a {\it uniform} approach, where
the same number of modes is adopted everywhere~\cite{PerottoZilio13,Perotto14b}. 
This flexibility in the choice of $m$ is currently not available for PGD.
Adaptive strategies to select the modal index are available for the three variants of the HiMod 
procedure~\cite{PerottoVeneziani14,PerottoZilio15}.

\subsection{Computational aspects}
\begin{figure}[t]
\hspace*{-.3cm}
\includegraphics[width=6.1cm,trim={2.9cm 0.5cm 2.7cm 0.5cm},clip]{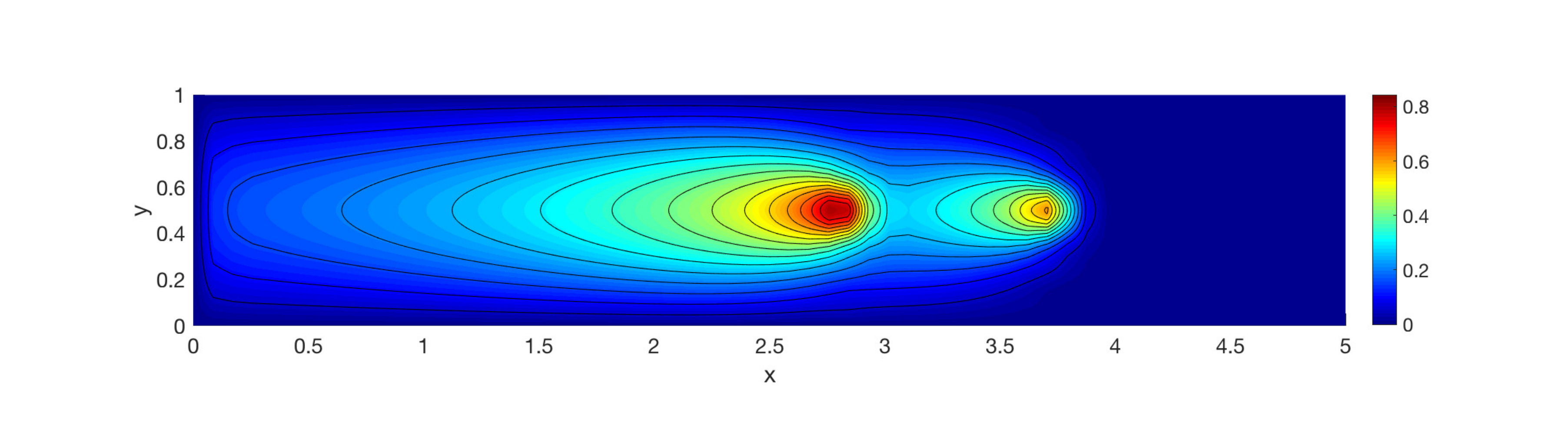}
\includegraphics[width=6.1cm,trim={3.2cm 0.5cm 2.4cm 0.5cm},clip]{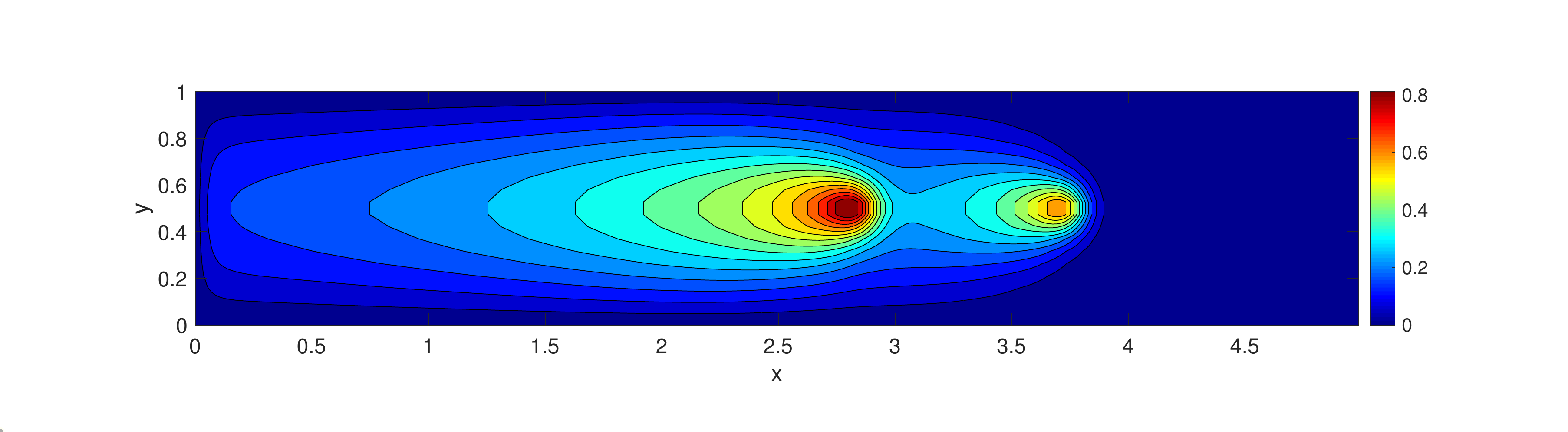}
\caption{Qualitative comparison between a HiMod (left) and a PGD (right) 
approximations.}\label{fig_HiModPGD}       
\end{figure}

From a computational viewpoint, HiMod reduction and PGD lead to completely different procedures.
Indeed, for a fixed value of $m$, we have to solve  
the only system \eqref{eq:himodalgebraic} of order $mN_h$ when applying HiMod,
in contrast to PGD which demands a multiple solution of systems \eqref{sys1}-\eqref{sys2} of order $N_h^x$ and $N_h^{\vec{y}}$, respectively
because of the fixed point and the enrichment algorithms.
Thus, the direct solution of a single system, in general of larger order, is replaced by an iterative solution of several and smaller systems.
This heterogeneity makes a computational comparison between PGD and HiMod not so meaningful.
We verify the reliability of the HiMod and PGD procedures on a common test case, by choosing in \eqref{eq:fullweak}
$V=H^1_0(\Omega)$ with $\Omega=(0, 5)\times(0, 1)$, 
$a(u,v)=\int_{\Omega} \big[ \mu \nabla u \cdot \nabla v + {\bf b} \cdot \nabla u \big]d\Omega$ 
for $\mu=0.24$, ${\bf b}=[-5, 0]^T$, and $F(v)=\int_{\Omega} fv d\Omega$ with 
$f(x,y) = 50 \big\{\exp\big[ - \big((x-2.85)/0.075 \big)^2 - \big( (y - 0.5)/0.075 \big)^2 \big] +
\exp\big[ - \big((x-3.75)/0.075 \big)^2 - \big( (y - 0.5)/0.075 \big)^2 \big] \big\}$.
For both the methods, we uniformly subdivide $\Omega_{1D}$ into $285$ subintervals.
We set the PGD discretization along $y$ as well as the PGD and the HiMod index $m$ in order to ensure the same accuracy, ${\tt TOL}$,
on the reduced approximations with respect to a reference FE solution, computed on a $2500\times 500$ structured mesh.
In particular, for ${\tt TOL}=8\cdot 10^{-3}$, we have to subdivide interval $(0,1)$ into $20$ uniform subintervals, and to set $m$ 
to $6$ and to $9$ in the PGD and the HiMod discretization, respectively. Sinusoidal functions are chosen for the HiMod modal basis.
The ADS iterations are controlled in terms of the relative increment, as
\begin{equation}\label{ads_cr}
\frac{\| u_m^{x, p} u_m^{{y}, p} - u_m^{x, p-1} u_m^{{y}, p-1} \|_{L^2(\Omega)}}{\| u_m^{x, p} u_m^{{y}, p} \|_{L^2(\Omega)}}
\le {\tt TOL_{FP}},
\end{equation}
with ${\tt TOL_{FP}}=10^{-2}$. Fig.~\ref{fig_HiModPGD} shows the reduced approximations (which are fully comparable with the FE one, here omitted).
The contourplots are very similar. The coarse PGD $y$-discretization justifies the slight
roughness of the PGD contourlines.

Another distinguishing feature between HiMod and PGD is the domain discretization. 
Indeed, HiMod requires only the partition $\mathcal{T}_h$ along $\Omega_{1D}$, independently of the dimension of $\Omega$. 
No discretization is needed in the $\vec{y}$-direction, although we have to carefully select the quadrature nodes 
to compute coefficients $\{ \hat r^{a, b}_{jk} \}$. This task becomes particularly challenging when dealing with polar coordinates~\cite{GuzzettiPerottoVeneziani18}.
With PGD to benefit of the computational advantages associated with a 1D discretization, we are obliged to assume the full separability
of $\Omega$; actually, a partial separability demands a 1D partition for $\Omega_x$, and a two-dimensional 
partition of $\Omega_{\vec{y}}$. As explained in Sect.~\ref{parametri_sec}, non-Cartesian domains require a 3D discretization of $\Omega$. 

Finally we analyze the interplay between the enrichment and the ADS iterations in the  
PGD reduction. 
We investigate the possible relationship between ${\tt TOL_{FP}}$ in \eqref{ads_cr} and ${\tt TOL_{E}}$ in \eqref{m_criterion}, to verify 
if a small tolerance for the fixed point iteration improves the accuracy 
of the PGD approximation, thus reducing the number of enrichment steps. To do this, 
we adopt the same test case used above. 
Table~\ref{table1} gathers the number of ADS iterations, \#$\tt IT_{FP}$, the number, $m$,
of enrichment steps, and the CPU time\footnote{The computations have 
been run on a Intel Core i5 Dual-Core CPU 2.7 GHz 8GB RAM MacBook.} (in seconds)
demanded by the PGD procedure, 
for two different values of ${\tt TOL_{E}}$ and three different choices of $\tt TOL_{FP}$. In particular, in column \#$\tt IT_{FP}$
we specify the number of ADS iterations required by each enrichment step. As expected, there exists a link between the two tolerances, namely,
when a higher accuracy constrains the fixed point iteration, a smaller number of enrichment steps is performed to ensure the accuracy
${\tt TOL_{E}}$.
\begin{table}[]
\caption{Quantitative analysis for PGD in terms of fixed point iterations and enrichment steps}\label{table1}       
\begin{tabular}{r|rlr|rlrlr|r|}
                                    & \multicolumn{3}{c|}{${\tt TOL_{E}}=2 \cdot 10^{-2}$} & \multicolumn{3}{c}{$ {\tt TOL_{E}}=8 \cdot 10^{-3}$} \\ \cline{2-7} 
                                    & \#$\tt IT_{FP}$   & $\quad m\quad $  & CPU {[}s{]}  & \#$\tt IT_{FP}$   & $\quad m\quad $  & CPU {[}s{]}  \\ \svhline
${\tt TOL_{FP} }= 10^{-1}$ & $\{ 2, 2, 2 \}$       &$\quad 3\quad$      & 0.099640     & $\{ 2,2,2,2,2 \}$     &$\quad 5\quad $      & 0.337861     \\
${\tt TOL_{FP} }= 10^{-2}$ & $\{ 4, 3 \}$           &$\quad 2\quad$      & 0.046756      & $\{ 4,3,2,2,4 \}$       & $\quad 5\quad $      & 0.358555     \\
${\tt TOL_{FP} }= 10^{-3}$ & $\{ 5, 5 \}$           &$\quad 2\quad$      & 0.077958     & $\{ 5,5,2,7 \}$          & $\quad 4\quad $      & 0.341748    
\end{tabular}
\end{table}

\section{HiMod reduction and PGD for parametrized problems}\label{parametri_sec}
The actual potential of PGD becomes more evident when considering a parametric setting, i.e., when problem \eqref{eq:fullweak} is replaced by the formulation
\begin{equation}\label{eq:fullweakparam}
	\text{find } u({\boldsymbol \mu}) \in V: a(u({\boldsymbol \mu}),v; {\boldsymbol \mu}) = F(v; {\boldsymbol \mu}) \quad \forall v \in V,
\end{equation}
with ${\boldsymbol \mu}$ a parameter, which may represent any data of the problem, 
e.g., the coefficients of the considered PDE, the source term,
a boundary value or the domain geometry.

The technique adopted by PGD to deal with the parametric dependence in \eqref{eq:fullweakparam}
is very effective. Parameter $\boldsymbol \mu$ is considered as an additional independent variable which
varies in a domain $\Omega_{\boldsymbol \mu}$~\cite{ChinestaKeuningsLeygue13}.
Thus, the PGD space \eqref{eq:pgdspace} changes into the new one
\begin{equation}\label{eq:pgdspace_para}
\begin{array}{lll}
W_m^{\boldsymbol \mu} = \Bigg\{ &w_m(x, \vec{y}, \boldsymbol \mu) = \displaystyle \sum_{k = 1}^m w_k^x(x) w_k^{\vec{y}}(\vec{y})
w_k^{\boldsymbol \mu}(\boldsymbol \mu), \text{ with }\\[2mm]  
& w_k^x \in W^x_h, \, w_k^{\vec{y}} \in W^{\vec{y}}_h,\, w_k^{\boldsymbol \mu} \in W^{\boldsymbol \mu}_h,  
x \in \Omega_x, \, \vec{y} \in \Omega_{\vec{y}}, \, \boldsymbol{\mu} \in \Omega_{\boldsymbol \mu} \Bigg\},
\end{array}
\end{equation}
with $W^{\boldsymbol \mu}_h$ a discretization of the space $L^2(\Omega_{\boldsymbol \mu}; \mathbb{R}^Q)$,
being $Q$ the length of vector $\boldsymbol \mu$.
Generalizing the enrichment paradigm in \eqref{pgd_rep}, at the $m$-th step of the PGD approach applied 
to problem \eqref{eq:fullweakparam} we have to compute 
three unknown functions, $u_m^x$, $u_m^{\vec{y}}$ and $u_m^{\boldsymbol \mu}$, by picking the test function as
$X(x)Y(\vec{y})Z({\boldsymbol \mu})$,
with $X\in W_h^x$, $Y\in W^{\vec{y}}_h$, $Z\in W_h^{\boldsymbol \mu}$. 
Functions $u_m^x$, $u_m^{\vec{y}}$, $u_m^{\boldsymbol \mu}$ are computed by ADS, which now coincides with a three-step procedure. Thus, with reference to the Poisson problem, $-\nabla \cdot \big( \mu \nabla u\big)=f$ completed with full homogeneous Dirichlet boundary conditions and for ${\boldsymbol \mu}\equiv\mu$, we first compute $u_m^{x, p}$ by identifying $u_m^{\vec{y}}$ and $u_m^{\boldsymbol \mu}$ with the previous approximations, $u^{\vec{y},p-1}_m$ and $u_m^{\boldsymbol \mu, p-1}$, respectively  
and by selecting $Y(\vec{y})Z({\boldsymbol \mu})=u^{\vec{y},p-1}_m u^{{\boldsymbol \mu},p-1}_m$ in the test function. This leads to a linear system which generalizes \eqref{sys1}, namely
\begin{equation}\label{sys1p}
\begin{array}{ll}
&\Big[ \big( {\bf u}_m^{{\boldsymbol \mu}, p-1}\big)^T M^{\boldsymbol \mu} {\bf u}_m^{{\boldsymbol \mu}, p-1} \Big] \Big\{ \Big[ \big( {\bf u}_m^{\vec{y}, p-1}\big)^T M^{\vec{y}} {\bf u}_m^{\vec{y}, p-1}\Big] K^{x} + \Big[ \big( {\bf u}_m^{\vec{y}, p-1}\big)^T K^{\vec{y}} {\bf u}_m^{\vec{y}, p-1}\Big] M^x \Big\} {\bf u}_m^{x, p} \\[2mm] 
= & \Big[ \big( {\bf u}_m^{\vec{y}, p-1}\big)^T {\bf f}^{\vec{y}}\Big] \Big[ \big( {\bf u}_m^{{\boldsymbol \mu}, p-1}\big)^T {\bf f}^{{\boldsymbol \mu}}\Big]  {\bf f}^{x}
- \sum_{k=1}^{m-1} \Big[ \big( {\bf u}_m^{{\boldsymbol \mu}, p-1}\big)^T M^{\boldsymbol \mu} {\bf u}_k^{{\boldsymbol \mu}} \Big]\\[2mm]
& \Big\{ \Big[ \big( {\bf u}_m^{\vec{y}, p-1}\big)^T M^{\vec{y}} {\bf u}_k^{\vec{y}}\Big] K^{x} + \Big[ \big( {\bf u}_m^{\vec{y}, p-1}\big)^T K^{\vec{y}} {\bf u}_k^{\vec{y}}\Big] M^x \Big\} {\bf u}_k^{x},
\end{array}
\end{equation}
where $M^{\boldsymbol \mu}\in \mathbb{R}^{N_h^{\boldsymbol \mu} \times N_h^{\boldsymbol \mu}}$ is the mass matrix associated with the parameter ${\boldsymbol \mu}$, with 
$\big[ M^{\boldsymbol \mu} \big]_{ij}=\int_{\Omega_{\boldsymbol \mu}} \boldsymbol \mu \theta_i^{\boldsymbol \mu} \theta_j^{\boldsymbol \mu} d{\boldsymbol \mu}$ for $i$, $j=1,\ldots, N_h^{\boldsymbol \mu}$ and 
${\mathcal B}_{\boldsymbol \mu}=\{ \theta_\gamma^{\boldsymbol \mu} \}_{\gamma = 1}^{N_h^{\boldsymbol \mu}}$ a basis for the space $W_h^{\boldsymbol \mu}$, ${\bf f}^{\boldsymbol \mu} \in \mathbb{R}^{N_h^{\boldsymbol \mu} }$ 
with $\big[ {\bf f}^{\boldsymbol \mu}\big]_l=\int_{\Omega_{\boldsymbol \mu}} f^{\boldsymbol \mu} \theta_l^{\boldsymbol \mu} d{\boldsymbol \mu}$ for $l=1,\ldots, N_h^{\boldsymbol \mu}$ after assuming the separability $f=f^x f^{\vec{y}} f^{\boldsymbol \mu}$
for the source term $f$, and where
we employ the same notation as in \eqref{sys1}-\eqref{sys2} 
to denote vectors ${\bf u}_w^{\boldsymbol \mu}$, ${\bf u}_m^{\boldsymbol \mu, s}$, with $w=1, \ldots, m-1$, $s=p, p-1$,
collecting the PGD coefficients associated with the basis ${\mathcal B}_{\boldsymbol \mu}$.
Analogously, $u_m^{\vec{y}, p}$ is computed by solving the generalization of the linear system \eqref{sys2} given by
$$
\begin{array}{ll}
&\Big[ \big( {\bf u}_m^{{\boldsymbol \mu}, p-1}\big)^T M^{\boldsymbol \mu} {\bf u}_m^{{\boldsymbol \mu}, p-1} \Big] \Big\{ \Big[ \big( {\bf u}_m^{x, p}\big)^T K^{x} {\bf u}_m^{x, p}\Big] M^{\vec{y}} + \Big[ \big( {\bf u}_m^{x, p}\big)^T M^{x} {\bf u}_m^{x, p}\Big] K^{\vec{y}} \Big\} {\bf u}_m^{\vec{y}, p} \\[2mm] 
= & \Big[ \big( {\bf u}_m^{x, p}\big)^T {\bf f}^{x}\Big] \Big[ \big( {\bf u}_m^{{\boldsymbol \mu}, p-1}\big)^T {\bf f}^{{\boldsymbol \mu}}\Big]  {\bf f}^{\vec{y}}
- \sum_{k=1}^{m-1} \Big[ \big( {\bf u}_m^{{\boldsymbol \mu}, p-1}\big)^T M^{\boldsymbol \mu} {\bf u}_k^{{\boldsymbol \mu}} \Big]\\[2mm]
& \Big\{ \Big[ \big( {\bf u}_m^{x, p}\big)^T K^{x} {\bf u}_k^{x}\Big] M^{\vec{y}} + \Big[ \big( {\bf u}_m^{x, p}\big)^T M^{x} {\bf u}_k^{x}\Big] K^{\vec{y}} \Big\} {\bf u}_k^{\vec{y}},
\end{array}
$$
after setting $u_m^x=u_m^{x, p}$, $u_m^{\boldsymbol \mu}=u_m^{\boldsymbol \mu, p-1}$ and 
$X(x)Z({\boldsymbol \mu})=u^{x,p}_m u^{{\boldsymbol \mu},p-1}_m$ for the PGD test function. 
Finally, we have the additional linear system used to compute $u_m^{\boldsymbol \mu, p}$,
$$
\begin{array}{ll}
& \Big\{ \Big[ \big( {\bf u}_m^{x, p}\big)^T K^{x} {\bf u}_m^{x, p}\Big] \Big[ \big( {\bf u}_m^{\vec{y}, p}\big)^T M^{\vec{y}} {\bf u}_m^{\vec{y}, p}\Big] + \Big[ \big( {\bf u}_m^{x, p}\big)^T M^{x} {\bf u}_m^{x, p}\Big] \Big[ \big( {\bf u}_m^{\vec{y}, p}\big)^T K^{\vec{y}} {\bf u}_m^{\vec{y}, p}\Big]
\Big\}  \\[2mm] 
& M^{\boldsymbol \mu}  {\bf u}_m^{\boldsymbol \mu, p}=  \Big[ \big( {\bf u}_m^{x, p}\big)^T {\bf f}^{x}\Big] \Big[ \big( {\bf u}_m^{\vec{y}, p}\big)^T {\bf f}^{\vec{y}}\Big]  {\bf f}^{\boldsymbol \mu}
- \sum_{k=1}^{m-1} \Big\{\Big[ \big( {\bf u}_m^{x, p}\big)^T K^{x} {\bf u}_k^{x} \Big]\Big[ \big( {\bf u}_m^{\vec{y}, p}\big)^T M^{\vec{y}} {\bf u}_k^{x} \Big]\\[2mm]
& \Big[ \big( {\bf u}_m^{x, p}\big)^T M^{x} {\bf u}_k^{x}\Big] \Big[ \big( {\bf u}_m^{\vec{y}, p}\big)^T K^{\vec{y}} {\bf u}_k^{\vec{y}}\Big]  \Big\} M^{\boldsymbol \mu} {\bf u}_k^{\boldsymbol \mu},
\end{array}
$$
obtained for $u_m^x=u_m^{x, p}$, $u_m^{\vec{y}}=u_m^{\vec{y}, p}$ and by selecting 
$X(x)Y(\vec{y})=u^{x,p}_m u^{\vec{y},p}_m$ for the test function.
From a computational viewpoint, at each ADS iteration, we have to solve now three linear systems 
of order $N_h^x$, $N_h^{\vec{y}}$, $N_h^{\boldsymbol \mu}$, respectively.

\begin{figure}[t]
\includegraphics[width=4.8cm]{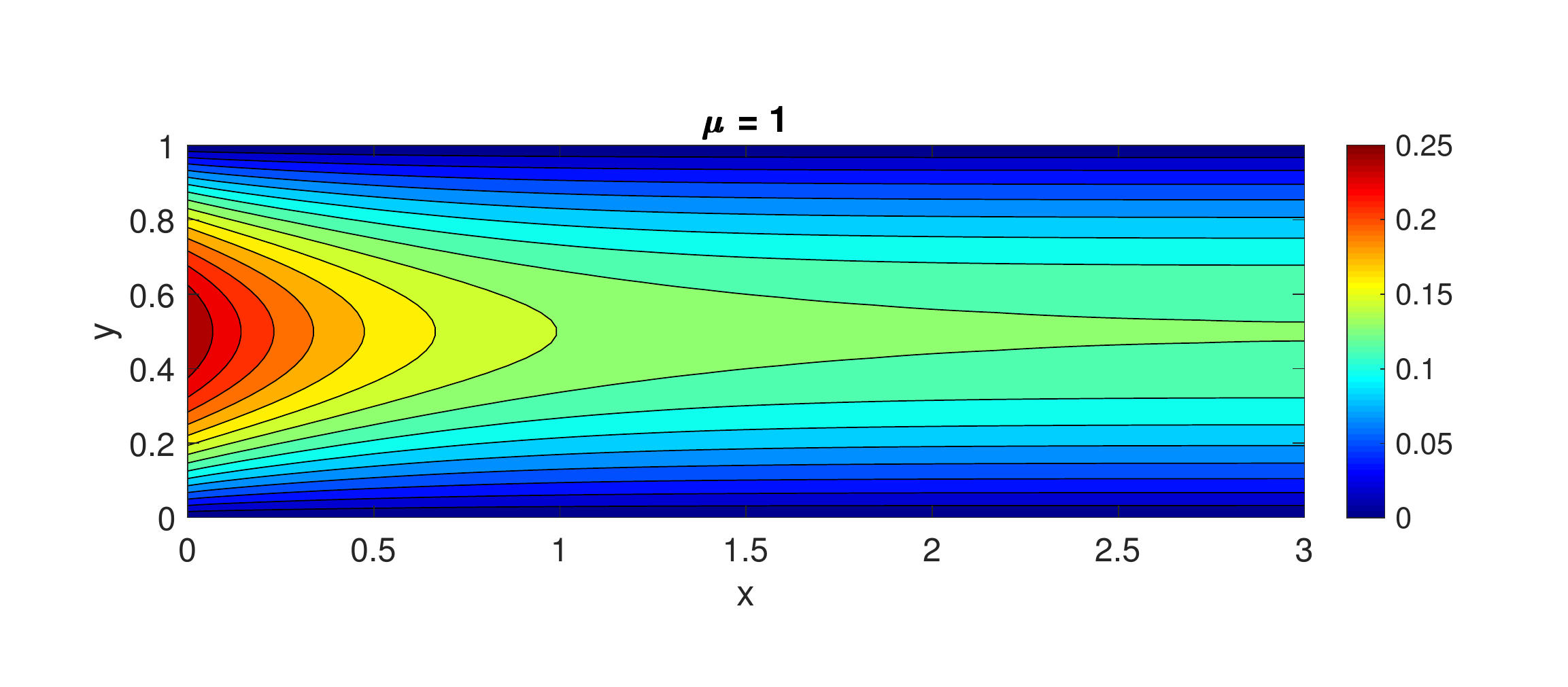}\hspace*{-.3cm}
\includegraphics[width=4.8cm]{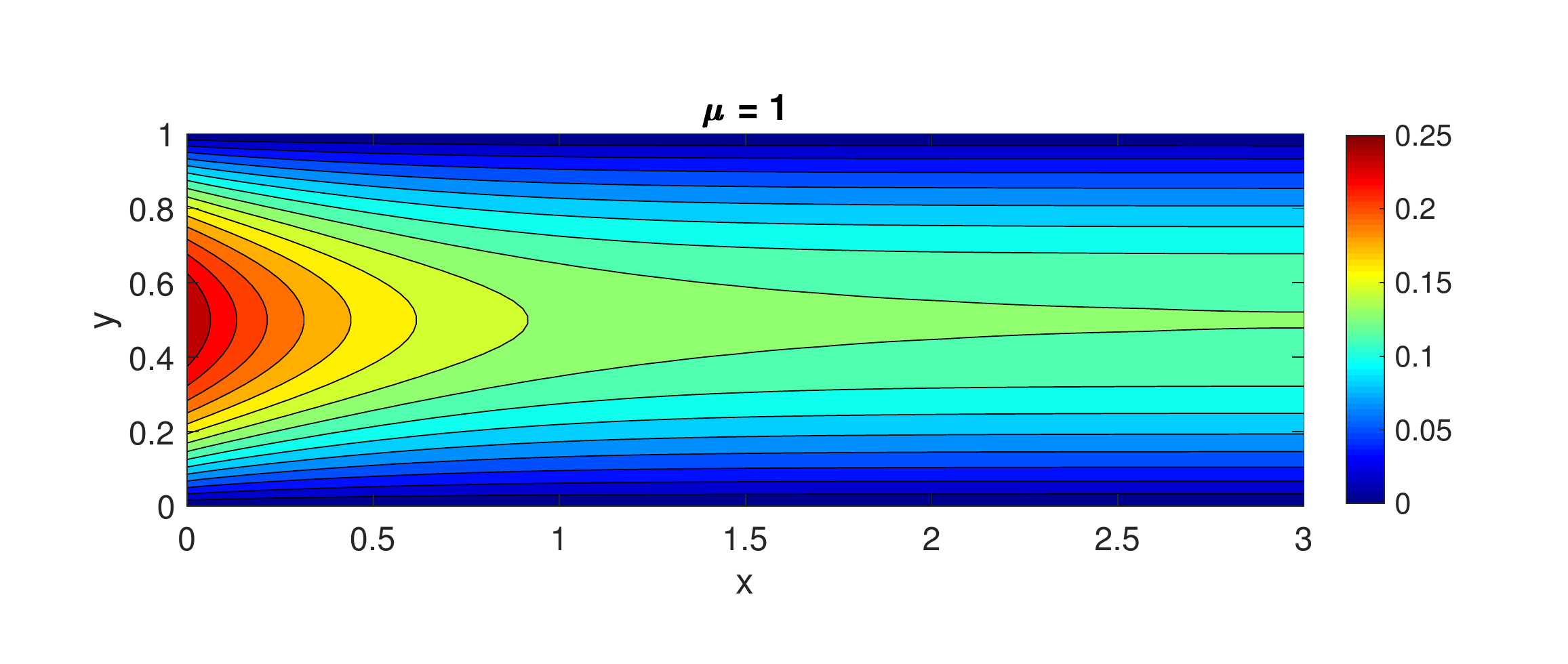}\\
\includegraphics[width=4.8cm]{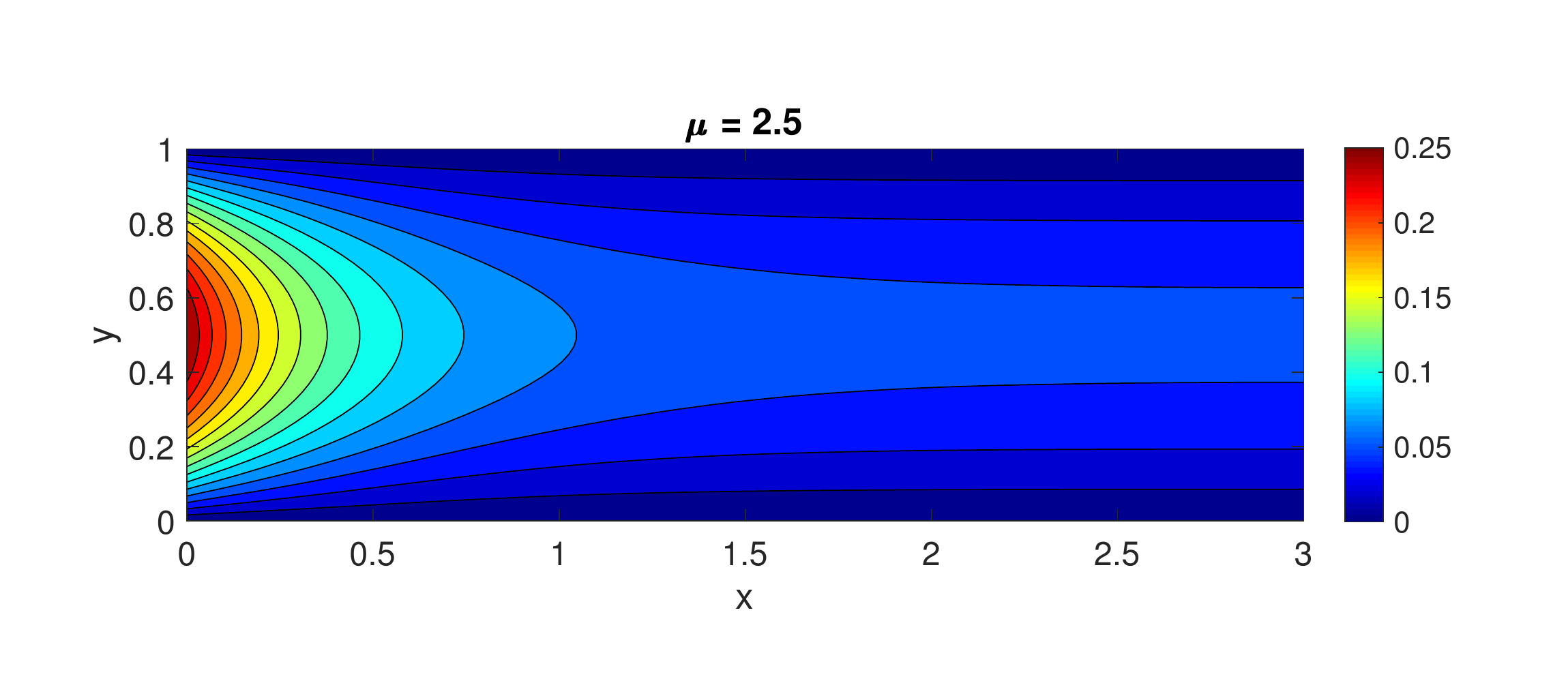}\hspace*{-.3cm}
\includegraphics[width=4.8cm]{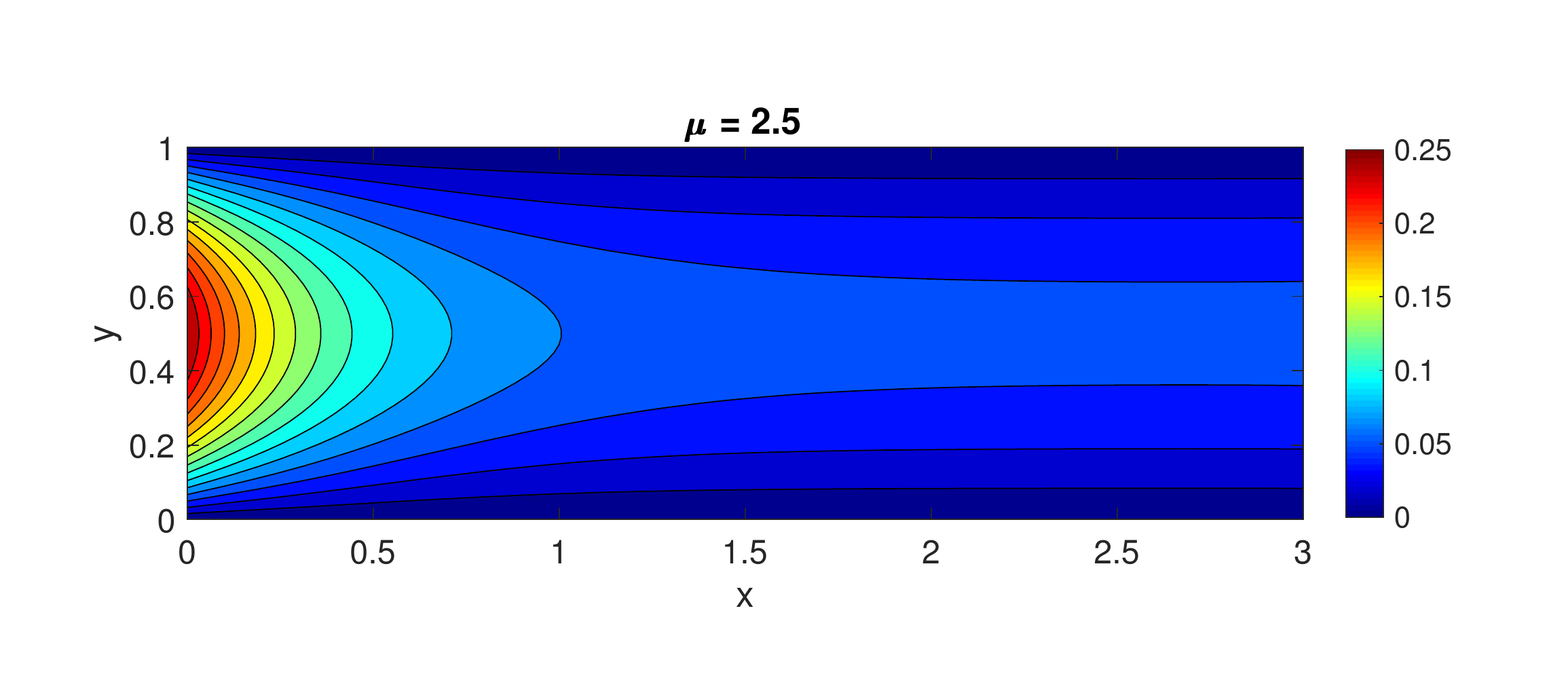}
\caption{Qualitative comparison between  the reference (left) and the PGD (right) 
solutions, for $\mu=1$ (top) and $\mu=2.5$ (bottom).}\label{fig_PGD}       
\end{figure}

We investigate the reliability of PGD on problem \eqref{eq:fullweakparam}, for 
$V=H^1_{\Gamma_{\rm in}\cup \Gamma_{\rm up}\cup \Gamma_{\rm down}}(\Omega)$
with $\Omega=(0,3) \times (0,1)$, $\Gamma_{\rm in}= \{0 \} \times (0,1)$, 
$\Gamma_{\rm up}=(0,3) \times \{ 1\}$, $\Gamma_{\rm down}=(0,3) \times \{ 0\}$,
$a(u,v)=\int_{\Omega} \big[ \mu \nabla u \cdot \nabla v + {\bf b} \cdot \nabla u \big]d\Omega$
with ${\bf b}=[2.5, 0]^T$ and $\mu$ the parameter to be varied in $\Omega_{\boldsymbol \mu}=[1, 5]$,
$F(v)=\int_{\Omega} fv d\Omega$ with $f=1$. The problem is completed with mixed boundary conditions, namely
a homogeneous Dirichlet data on $\Gamma_{\rm up}\cup \Gamma_{\rm down}$, the non-homogeneous 
Dirichlet condition, $u=u_{\rm in}$ with $u_{\rm in}=y(1-y)$, on $\Gamma_{\rm in}$ and a homogeneous Neumann value on
$\Gamma_{\rm out}=\{ 3 \} \times (0,1)$.
We apply the PGD reduction for $m=2$, and we uniformly subdivide $\Omega_x$, $\Omega_{y}$,
$\Omega_{\mu}$, being $N_h^x=150$, $N_h^{y}=50$, $N_h^{\mu}=500$.
The tolerance in \eqref{ads_cr} is set to $10^{-2}$. 
Fig.~\ref{fig_PGD} compares the PGD approximation for $\mu=1$ and $\mu=2.5$
with a reference full solution coinciding with a linear FE approximation computed on a $300\times 100$
structured mesh. The qualitative matching between the corresponding solutions is significant.
From a quantitative viewpoint, the $L^2(\Omega)$-norm of the relative error associated with the PGD approximation does not vary significantly
by increasing $m$, whereas a slight error reduction is detected by increasing $\mu$. 

The parametric counterpart of the HiMod reduction, known as HiPOD,  
merges HiMod with POD~\cite{BarolietAl17,LupoPasiniPerottoVeneziani19}. HiPOD pursues a different goal with 
respect to PGD. Indeed, for a new value, ${\boldsymbol \mu}^*$, of the parameter, PGD provides an approximation for the 
full solution $u({{\boldsymbol \mu}^*}) $, while HiPOD approximates the HiMod solution associated with ${\boldsymbol \mu}^*$.
The offline/online paradigm of POD is followed also by HiPOD. 
The peculiarity is that the offline step is now performed in the HiMod setting
to contain the computational burden typical of this stage and by relying on 
the good properties of HiMod in terms of reliability-versus-accuracy balance.
Thus, we choose $P$ different values, 
$\boldsymbol{\mu}=\boldsymbol{\mu}_i$ with $i=1, \ldots, P$,
for parameter $\boldsymbol{\mu}$, and we collect the HiMod approximation 
for the corresponding problem \eqref{eq:fullweakparam} into the response matrix,
$\mathcal S = \big[  \vec{u}_m^{\rm HiMod} (\boldsymbol{\mu}_1), \vec{u}_m^{\rm HiMod} (\boldsymbol{\mu}_2), \ldots, 
\vec{u}_m^{\rm HiMod} (\boldsymbol{\mu}_P) \big]\in \mathbb{R}^{mN_h\times P}$, according to representation \eqref{mod_expansion}. Successively, we define the null-average matrix 
$$
\mathcal V=\mathcal S - \frac{1}{P}\sum_{i=1}^P \left[ \vec{u}_m^{\rm HiMod} (\boldsymbol{\mu}_i), \vec{u}_m^{\rm HiMod} (\boldsymbol{\mu}_i), \ldots, \vec{u}_m^{\rm HiMod} (\boldsymbol{\mu}_i) \right]\in \mathbb{R}^{mN_h\times P},
$$
and we apply the Singular Value Decomposition (SVD) to $\mathcal V$, so that $\mathcal V=\Phi \Sigma \Psi^T$,
where $\Phi\in \mathbb{R}^{(mN_h) \times (mN_h)}$ and $\Psi\in \mathbb{R}^{P\times P}$ are the unitary matrices of the 
left- and of the right-singular vectors of ${\mathcal V}$, respectively while
$\Sigma = \text{ diag }(\sigma_1, \dots , \sigma_{\rho}) \in \mathbb{R}^{(mN_h) \times P}$ denotes the pseudo-diagonal matrix
of the singular values of ${\mathcal V}$, being  
$\sigma_1 \geq \sigma_2 \geq \dots \geq \sigma_{\rho}\ge 0$ and $\rho = \min(m N_h, P)$~\cite{GolubVanLoan13}.
The POD basis is identified by the first $l$ left singular 
vectors, ${\boldsymbol \phi}_i$, of $\mathcal V$, so that the reduced POD space is
$V_{\rm POD}^l= \text{span} \{{\boldsymbol \phi}_1, \dots, {\boldsymbol \phi}_{l}\}$, with $\dim(V_{\rm POD}^l)=l$
and $l\ll mN_h$. In the numerical assessment below, 
value $l$ coincides with the smallest integer such that $\sigma_{l}^2< \varepsilon$, with $\varepsilon$ 
a prescribed tolerance. \\
The online phase of HiPOD approximates the HiMod solution to problem \eqref{eq:fullweakparam}
for a new value, ${\boldsymbol{\mu}}^*$, of the parameter by exploiting the POD basis instead of solving system \eqref{eq:himodalgebraic}. 
This is performed via a projection step. After assembling the HiMod stiffness matrix and right-hand side,
$A_m^{\rm HiMod}({\boldsymbol{\mu}}^*)$ and $\vec{f}_m^{\rm HiMod}({\boldsymbol{\mu}}^*)$, associated with the new value
of the parameter, we solve the POD system of order $l$
\begin{equation}\label{pod_s}
A_{\rm POD}({\boldsymbol{\mu}}^*) {\bf u}_{\rm POD}({\boldsymbol{\mu}}^*)={\bf f}_{\rm POD}({\boldsymbol{\mu}}^*),
\end{equation}
where $A_{\rm POD}({\boldsymbol{\mu}}^*)=(\Phi_{\rm POD}^l)^T A_m^{\rm HiMod}({\boldsymbol{\mu}}^*) \, \Phi_{\rm POD}^l$ and ${\bf f}_{\rm POD}({\boldsymbol{\mu}}^*) = (\Phi_{\rm POD}^l)^T \mathbf{f}_m^{\rm HiMod}({\boldsymbol{\mu}}^*)$ denote the POD stiffness matrix and right-hand side, respectively with $\Phi_{\rm POD}^l=[{\boldsymbol \phi}_1, \ldots, {\boldsymbol \phi}_l]\in \mathbb{R}^{(mN_h) \times l}$ the matrix 
collecting the POD basis vectors. The HiMod solution is thus approximated by vector $\Phi_{\rm POD}^l {\bf u}_{\rm POD}({\boldsymbol{\mu}}^*)\in \mathbb{R}^{mN_h}$, i.e., after solving a system of order $l$ instead of $mN_h$.
Overall, HiPOD requires to solve $P$ linear systems of order $mN_h$ during the offline phase,
additionally to a system of order $l$ in the online phase.
\begin{figure}[t]
  \begin{minipage}[t]{0.39\textwidth}
    \includegraphics[width=\textwidth]{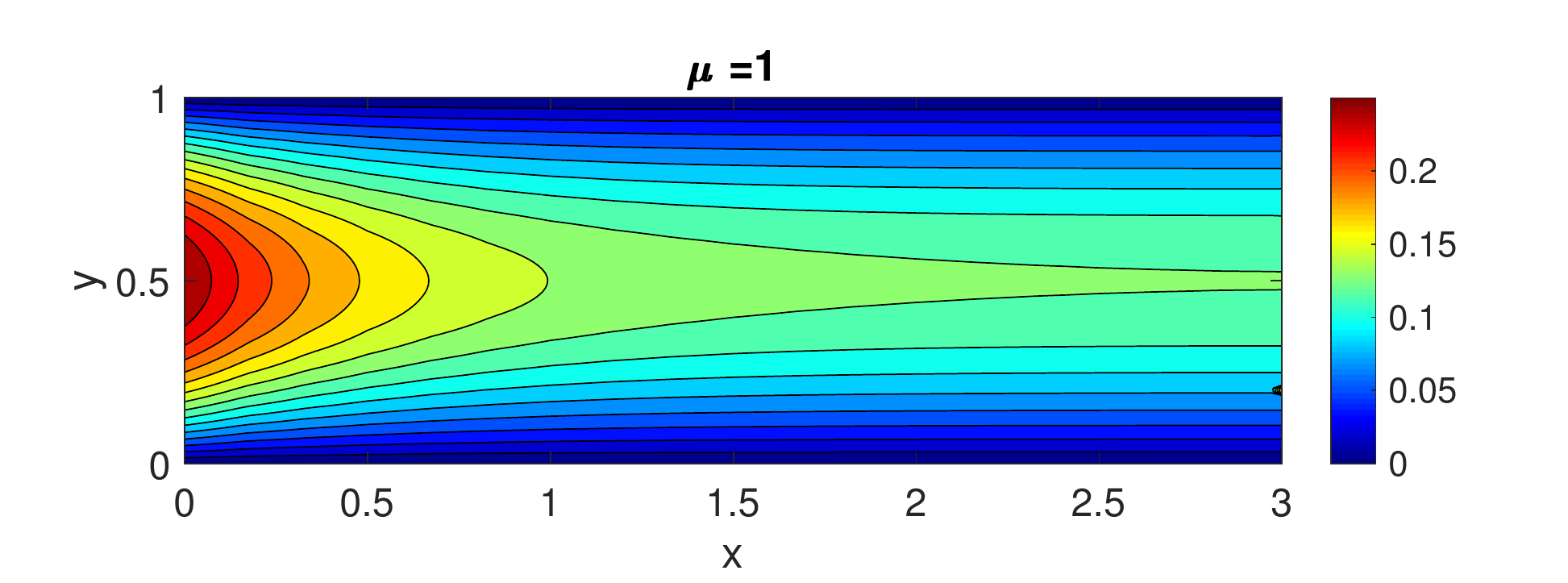}\\
    \includegraphics[width=\textwidth]{{{himod_25}}}
  \end{minipage}
  \begin{minipage}[t]{0.39\textwidth}
    \hspace*{-.45cm}\includegraphics[width=\textwidth]{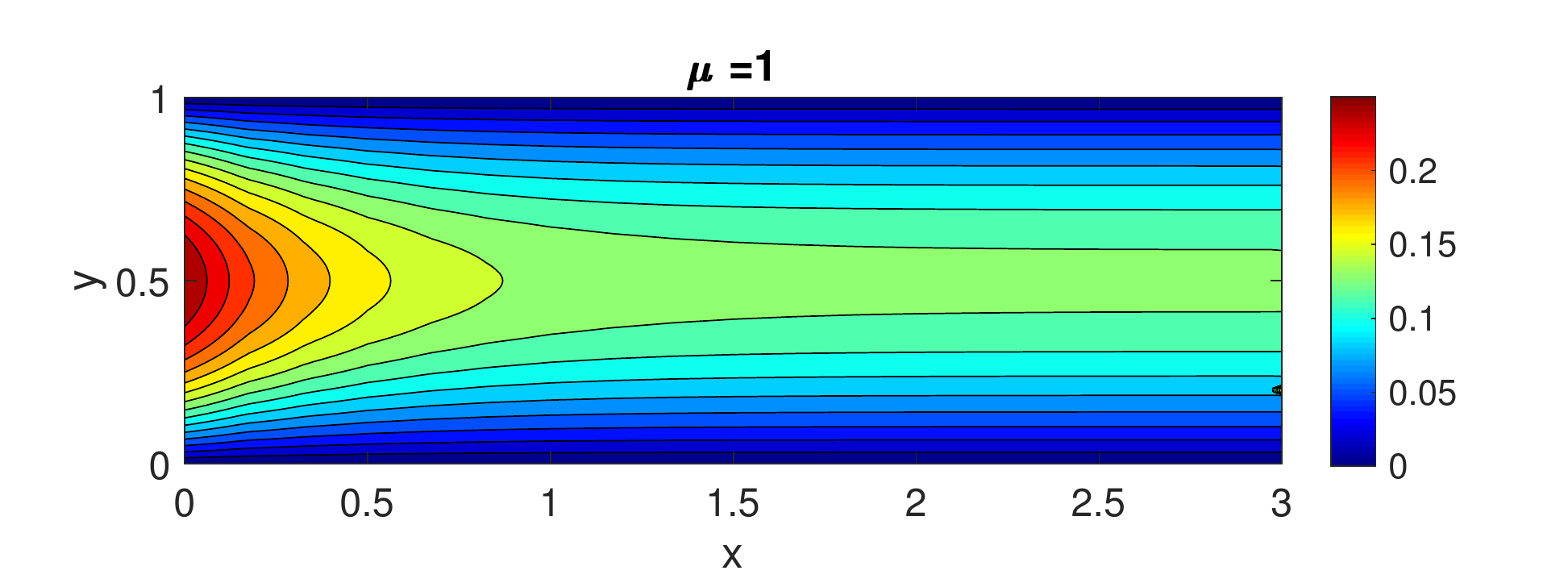}\\
    \hspace*{-.45cm}\includegraphics[width=\textwidth]{{{hipod_mu_25}}}
  \end{minipage}
  \begin{minipage}[t]{0.1\textwidth}
   \begin{scriptsize}
   \hspace*{-.6cm}
   \begin{tabular}{lll}
                & $l=1$    & $l=4$   \\\hline
    $\mu^* = 1$   & 4.06e-02 & 2.53e-06 \\\hline
    $\mu^* = 2.5$ & 2.74e-03 & 1.11e-07 \\\hline
    random      & 7.19e-03 & 2.65e-07 \\\hline\\[2mm]\hline
                & $l=6$    & $l=8$   \\\hline
    $\mu^* = 1$   & 1.79e-09 & 5.58e-12 \\ \hline
    $\mu^* = 2.5$ & 4.05e-10 & 1.21e-13 \\ \hline
    random      & 2.97e-10 & 3.66e-13 \\ \hline
    \end{tabular}
    \end{scriptsize}
  \end{minipage}
\caption{Contour plots: comparison between  the reference HiMod solution (left) and the HiPOD approximation with $l=1$ (right), 
for $\mu^*=1$ (top) and $\mu*=2.5$ (bottom). Table:
relative error between HiMod and HiPOD solutions with respect to the $L^2(\Omega)$-norm.}\label{fig_HiPOD} 
\end{figure}

To check the performances of HiPOD, we adopt the test case used above for PGD, 
for the same values of the parameters, $\mu^*=1$ and $\mu^*=2.5$. The reference solution is 
the corresponding HiMod approximation computed by using $m=15$ sinusoidal functions
in the $y$-direction, and a linear FE discretization along the mainstream based on a uniform subdivision of $\Omega_{1D}$ into $50$ subintervals. The same HiMod discretization is adopted to build the response matrix.
Concerning the HiPOD approximation, we pick $P=100$ by uniformly sampling the interval $[1, 5]$, and we select 
$\varepsilon=2.5 \cdot 10^{-15}$.  This choice sets the dimension of the POD space to $l=8$, so that we have to solve a system
of order $8$ instead of $750$. The contour plots in Fig.~\ref{fig_HiPOD} qualitatively compare the HiMod solution with the HiPOD approximation for $l=1$.  
The correspondence between the two approximations is good despite a single POD mode 
is employed (in such a case, system \eqref{pod_s} reduces to a scalar equation). We do not provide the HiPOD approximations for $l=8$ 
since they qualitatively coincide with the corresponding HiMod solution. 
The left panels can be additionally compared with the FE solutions in Fig.~\ref{fig_PGD}
to verify the reliability of the HiMod procedure. Finally, the table in Fig.~\ref{fig_HiPOD} gathers the $L^2(\Omega)$-norm 
of the relative error between HiMod and HiPOD solutions, for four different POD bases and for three choices of the viscosity
($1$, $2.5$ and the average over a sampling of $30$ random values of $\mu$). 
The error monotonically decreases for larger and larger values of $l$, independently of the choice for $\mu$.
If we compare the values for $\mu=1$ and for $\mu=2.5$ 
(one of the endpoints and the midpoint of the sampling interval, respectively), we notice a higher accuracy 
(of about one order of magnitude)
for the latter choice. This is rather standard in projection-based reduced order modeling~\cite{HesthavenRozzaStamm16}.
Concerning the computational saving in terms of CPU time,  
HiPOD method requires  on average $O(10^{-3})$[s] to be compared with $O(10)$[s] demanded by HiMod, 
resulting in a speedup of $10^4$.

Although PGD and HiPOD are not directly comparable due to the different purpose they pursue, we highlight the main pros and cons of the two methods. 
The explicit dependence of the approximation on the parameters makes PGD an ideal tool to efficiently deal with parametric problems.
For any new parameter, a direct evaluation yields the corresponding PGD approximation.
On the other hand, HiPOD suffers of the drawbacks typical of the projection-based methods. The main bottleneck is the 
assembling of the HiMod arrays involved in $A_{\rm POD}({\boldsymbol{\mu}}^*)$ and ${\bf f}_{\rm POD}({\boldsymbol{\mu}}^*)$.\\
When PGD is applied to parametric problems, we recover the possibility to deal with any geometric domain. 
In such a case, a partial separability is applied to the problem, so that the space independent variables are kept together 
whereas parameters are separated. This approach clearly looses the computational advantages due to
space separability. On the contrary, HiPOD inherits the geometric flexibility 
of the HiMod reduction, without giving up the spatial dimensional reduction of the problem.

\begin{acknowledgement}
The authors thank Yves Antonio Brandes Costa Barbosa for his support in the HiMod simulations. 
This work has been partially funded by
GNCS-INdAM 2018 project on ``Tecniche di Riduzione di Modello per le Applicazioni Mediche''.
F. Ballarin also acknowledges the support by European Union Funding for Research and Innovation, Horizon 2020 Program, in the framework of European Research Council Executive Agency: H2020 ERC Consolidator Grant 2015 AROMA-CFD project 681447 ``Advanced Reduced Order Methods with Applications in Computational Fluid Dynamics'' (P.I. G. Rozza).
\end{acknowledgement}

\bibliographystyle{plain}
\bibliography{referenc}

\begin{thebibliography}{10}

\bibitem{AlettiPerottoVeneziani18}
M.C. Aletti, S.~Perotto, and A.~Veneziani.
\newblock Hi{M}od reduction of advection-diffusion-reaction problems with
  general boundary conditions.
\newblock {\em J. Sci. Comput.}, 76(1):89--119, 2018.

\bibitem{AmmaretAl10}
A.~Ammar, F.~Chinesta, P.~Diez, and A.~Huerta.
\newblock An error estimator for separated representations of highly
  multidimensional models.
\newblock {\em Comput. Methods Appl. Mech. Engrg.}, 199(25-28):1872--1880,
  2010.

\bibitem{AmmarCuetoChinesta12}
A.~Ammar, E.~Cueto, and F.~Chinesta.
\newblock Reduction of the chemical master equation for gene regulatory
  networks using proper generalized decompositions.
\newblock {\em Int. J. Numer. Methods Biomed. Eng.}, 28(9):960--973, 2012.

\bibitem{BarolietAl17}
D.~Baroli, C.M. Cova, S.~Perotto, L.~Sala, and A.~Veneziani.
\newblock Hi-{P}{O}{D} solution of parametrized fluid dynamics problems:
  preliminary results.
\newblock In {\em Model Reduction of Parametrized Systems}, volume~17 of {\em
  MS\&A. Model. Simul. Appl.}, pages 235--254. Springer, Cham, 2017.

\bibitem{ChinestaKeuningsLeygue13}
F.~Chinesta, R.~Keunings, and A.~Leygue.
\newblock {\em The Proper Generalized Decomposition for Advanced Numerical
  Simulations: a Primer}.
\newblock SpringerBriefs in Applied Sciences and Technology. Springer
  International Publishing, 2014.

\bibitem{ErnPerottoVeneziani08}
A.~Ern, S.~Perotto, and A.~Veneziani.
\newblock Hierarchical model reduction for advection-diffusion-reaction
  problems.
\newblock In {\em Numerical Mathematics and Advanced Applications}, pages
  703--710. Springer, Berlin, 2008.

\bibitem{GhantiosetAl12}
C.~Ghnatios, A.~Ammar, A.~Cimetiere, A.~Hamdouni, A.~Leygue, and F.~Chinesta.
\newblock First steps in the space separated representation of models defined
  in complex domains.
\newblock In {\em 11th Biennial Conference on Engineering Systems Design and
  Analysis}, pages 37--42. Nantes, 2012.

\bibitem{GolubVanLoan13}
G.H. Golub and C.F. Van~Loan.
\newblock {\em Matrix Computations}.
\newblock Johns Hopkins Studies in the Mathematical Sciences. Johns Hopkins
  University Press, Baltimore, MD, fourth edition, 2013.

\bibitem{GonzalesetAl10}
D.~Gonz\'{a}lez, A.~Ammar, F.~Chinesta, and E.~Cueto.
\newblock Recent advances on the use of separated representations.
\newblock {\em Internat. J. Numer. Methods Engrg.}, 81(5):637--659, 2010.

\bibitem{GuzzettiPerottoVeneziani18}
S.~Guzzetti, S.~Perotto, and A.~Veneziani.
\newblock Hierarchical model reduction for incompressible fluids in pipes.
\newblock {\em Internat. J. Numer. Methods Engrg.}, 114(5):469--500, 2018.

\bibitem{HesthavenRozzaStamm16}
J.S. Hesthaven, G.~Rozza, and B.~Stamm.
\newblock {\em Certified Reduced Basis Methods for Parametrized Partial
  Differential Equations}.
\newblock SpringerBriefs in Mathematics. Springer, Cham; BCAM Basque Center for
  Applied Mathematics, Bilbao, 2016.

\bibitem{LadevezePassieuxNeron10}
P.~Ladev\`eze, J.-C. Passieux, and D.~N\'{e}ron.
\newblock The {LATIN} multiscale computational method and the proper
  generalized decomposition.
\newblock {\em Comput. Methods Appl. Mech. Engrg.}, 199(21-22):1287--1296,
  2010.

\bibitem{LupoPasiniPerottoVeneziani19}
M.~Lupo~Pasini, S.~Perotto, and A.~Veneziani.
\newblock {H}i{P}{O}{D}: {H}ierarchical model reduction driven by a {P}roper
  {O}rthogonal {D}ecomposition for parametrized advection-diffusion-reaction
  problems.
\newblock In preparation.

\bibitem{NiroomandietAl13}
S.~Niroomandi, D.~Gonz\'{a}lez, I.~Alfaro, F.~Bordeu, A.~Leygue, E.~Cueto, and
  F.~Chinesta.
\newblock Real-time simulation of biological soft tissues: a {PGD} approach.
\newblock {\em Int. J. Numer. Methods Biomed. Eng.}, 29(5):586--600, 2013.

\bibitem{Perotto14}
S.~Perotto.
\newblock Hierarchical model ({H}i-{M}od) reduction in non-rectilinear domains.
\newblock In {\em Domain Decomposition Methods in Science and Engineering
  {XXI}}, volume~98 of {\em Lect. Notes Comput. Sci. Eng.}, pages 477--485.
  Springer, Cham, 2014.

\bibitem{Perotto14b}
S.~Perotto.
\newblock A survey of hierarchical model ({H}i-{M}od) reduction methods for
  elliptic problems.
\newblock In {\em Numerical simulations of coupled problems in engineering},
  volume~33 of {\em Comput. Methods Appl. Sci.}, pages 217--241. Springer,
  Cham, 2014.

\bibitem{PerottoErnVeneziani10}
S.~Perotto, A.~Ern, and A.~Veneziani.
\newblock Hierarchical local model reduction for elliptic problems: a domain
  decomposition approach.
\newblock {\em Multiscale Model. Simul.}, 8(4):1102--1127, 2010.

\bibitem{PerottoRealiRusconiVeneziani17}
S.~Perotto, A.~Reali, P.~Rusconi, and A.~Veneziani.
\newblock H{IGAM}od: a hierarchical isogeometric approach for model reduction
  in curved pipes.
\newblock {\em Comput. \& Fluids}, 142:21--29, 2017.

\bibitem{PerottoVeneziani14}
S.~Perotto and A.~Veneziani.
\newblock Coupled model and grid adaptivity in hierarchical reduction of
  elliptic problems.
\newblock {\em J. Sci. Comput.}, 60(3):505--536, 2014.

\bibitem{PerottoZilio13}
S.~Perotto and A.~Zilio.
\newblock Hierarchical model reduction: three different approaches.
\newblock In {\em Numerical mathematics and advanced applications 2011}, pages
  851--859. Springer, Heidelberg, 2013.

\bibitem{PerottoZilio15}
S.~Perotto and A.~Zilio.
\newblock Space-time adaptive hierarchical model reduction for parabolic
  equations.
\newblock {\em Adv. Model. and Simul. in Eng. Sci.}, 2:25, 2015.

\bibitem{PruliereChinestaAmmar10}
E.~Pruliere, F.~Chinesta, and A.~Ammar.
\newblock On the deterministic solution of multidimensional parametric models
  using the proper generalized decomposition.
\newblock {\em Math. Comput. Simulation}, 81(4):791--810, 2010.

\bibitem{SignoriniZlotnikDiez17}
M.~Signorini, S.~Zlotnik, and P.~D\'{i}ez.
\newblock Proper generalized decomposition solution of the parameterized
  {H}elmholtz problem: application to inverse geophysical problems.
\newblock {\em Internat. J. Numer. Methods Engrg.}, 109(8):1085--1102, 2017.

\end{thebibliography}

\end{document}